\documentclass[11pt,leqno]{article}
\usepackage{amsmath,amsthm,txfonts}
\usepackage[ps2pdf,colorlinks=true,urlcolor=blue,citecolor=red,
linkcolor=blue,linktocpage,pdfpagelabels,bookmarksnumbered,
bookmarksopen]{hyperref}
\usepackage[english]{babel}

\numberwithin{equation}{section}

\newcommand{\eps}{\varepsilon}
\newcommand{\dt}{\delta}

\newcommand{\iu}{{\rm i}}

\newcommand{\la}{\lambda}
\newcommand{\bt}{\beta}

\newcommand{\vth}{\vartheta}

\renewcommand{\L}{\varmathbb{L}^{2}}
\newcommand{\E}{\mathcal{E}}
\newcommand{\I}{\mathcal{I}}
\renewcommand{\H}{\varmathbb{H}^{1}}
\newcommand{\M}{\mathcal{M}}

\newcommand{\ynu}{y_{m,1}}
\newcommand{\ynd}{y_{m,2}}
\newcommand{\unu}{u_{m,1}}
\newcommand{\und}{u_{m,2}}
\def\hsob{H^{1}(\R)}
\newcommand{\nehari}{\mathcal N}

\newcommand\N{\varmathbb{N}}
\newcommand\R{\varmathbb{R}}
\newcommand\C{\varmathbb{C}}
\newcommand\huno{H^{1}}
\newcommand\ldue{L^{2}}

\newcommand{\rdue}{\R^{2}}
\newcommand{\pa}{\partial}

\newtheorem{theorem}{Theorem}[section]
\newtheorem{remark}[theorem]{Remark}
\newtheorem{lemma}[theorem]{Lemma}
\newtheorem{proposition}[theorem]{Proposition}
\newtheorem{corollary}[theorem]{Corollary}
\newtheorem{definition}[theorem]{Definition}
\newcommand{\bos}{\begin{remark}\rm}
\newcommand{\eos}{\end{remark}}
\newcommand{\bte}{\begin{theorem}}
\newcommand{\ete}{\end{theorem}}
\newcommand{\bpr}{\begin{proposition}}
\newcommand{\epr}{\end{proposition}}
\newcommand{\bdf}{\begin{definition}\rm}
\newcommand{\edf}{\end{definition}}
\newcommand{\bco}{\begin{corollary}}
\newcommand{\eco}{\end{corollary}}
\newcommand{\ble}{\begin{lemma}}
\newcommand{\ele}{\end{lemma}}
\newcommand{\bdm}{\begin{displaymath}}
\newcommand{\edm}{\end{displaymath}}
\newcommand{\beq}{\begin{equation}}
\newcommand{\eeq}{\end{equation}}
\newcommand{\bdim}{{\noindent{\bf Proof.}\,\,\,}}
\newcommand{\edim}{{\unskip\nobreak\hfil\penalty50
 \hskip2em\hbox{}\nobreak\hfil\mbox{\rule{1ex}{1ex} \qquad}
  \parfillskip=0pt \finalhyphendemerits=0\par\medskip}}

\newcommand{\ov}{\overline}

\newcommand{\dys}{\displaystyle }

\begin{document}

\title{Energy convexity estimates for non-degenerate ground states of nonlinear 1D Schr\"odinger systems}

\author{ Eugenio Montefusco\thanks{Research supported by the
MIUR national research project\newline {\it Variational Methods and Nonlinear Differential Equations}.},\; Benedetta Pellacci$^{*}$,\;
Marco Squassina\thanks{Research supported by the MIUR national
research project\newline {\it``Variational and Topological Methods in the Study of
Nonlinear Phenomena''}.
 }}

\date{\today}
\maketitle

\begin{abstract}
We study the spectral structure of the complex linearized operator 
for a class of nonlinear Schr\"odinger systems, obtaining as  
byproduct some interesting properties of non-degenerate ground state
of the associated elliptic system, such as being isolated and 
orbitally stable.
\end{abstract}

\medskip

\section{Introduction and main results}

In the last few years, the interest in the study of Schr\"odinger systems
has considerably increased, in particular, for the following class 
of two weakly coupled nonlinear Schr\"odinger equations
\beq\label{schr}
\begin{cases}
\dys \iu\pa_{t}\phi_{1}+ \frac12\pa_{xx}\phi_{1}+\big(|\phi_{1}|^{2p}
+\beta|\phi_{2}|^{p+1}|\phi_{1}|^{p-1}\big)\phi_{1}=0 & \text{in $\R\times \R^{+}\!\!,$}
\smallskip
\\ 
\dys \iu\pa_{t}\phi_{2}+\frac12\pa_{xx}\phi_{2}+\big(|\phi_{2}|^{2p}
  +\beta|\phi_{1}|^{p-1}|\phi_{2}|^{p+1}\big)\phi_{2}=0 & 
\text{in $\R\times \R^{+}\!\!,$}
\\
\phi_{1}(0,x)=\phi_{1}^{0}(x), \quad \phi_{2}(0,x)=\phi_{2}^{0}(x)
& \text{in $\R$},  
\end{cases}
\eeq
where $\Phi=(\phi_{1},\phi_{2})$ and $\phi_{i}:[0,\infty)\times\R\to\C$, 
$\phi_{i}^{0}:\R\to\C$,  $0< p<2$. 
Usually the coupling constant $\beta>0$ models the birefringence 
effects inside a given anisotropic material (see e.g. \cite{manak}, \cite{men}). 
A soliton or standing wave solution is a solution of the form 
$\Phi(x,t)=(u_{1}(x)e^{\iu t},u_{2}(x)e^{\iu t})$ where 
$U(x)=(u_{1}(x),u_{2}(x))$ solves the elliptic system
\beq\label{ellittico}
\begin{cases}
-\dys \frac12  \pa_{xx} r_{1}+r_{1}=r_{1}^{2p+1}
  +\bt r_{1}^{p} r_{2}^{p+1}  & \text{in $\R$}, \smallskip\\
\dys  -\frac12\pa_{xx} r_{2}+r_{2}=r_{2}^{2p+1}
  +\bt r_{2}^{p} r_{1}^{p+1} & \text{in $\R$}.
\end{cases}
\eeq
Among all the solutions of \eqref{ellittico} there are the ground states, namely
least energy solutions. It is known  (see e.g. \cite{mmp1}, \cite{si}) that for $p\geq 1$ there exists a ground state $R=(r_{1},r_{2})$ $\in C^{2}(\R)\cap W^{2,s}(\R)$ for any positive $s$;  Moreover, $R$ has nonnegative components $r_i$ which are even,  decreasing on $\R^{+}$ and
exponentially decaying.
In \cite{mmp3} it is shown that $R$ can be characterized  as a  solutions of the following minimization problem 
\begin{equation}\label{elle2}
\E(R)=\inf_{\M}\E(V)\qquad \text{where }\quad 
\M:=\left\{V\in H^{1}(\R)\times H^{1}(\R),\,\|V\|_{2}=\|R\|_{2} \right\},
\end{equation}
and
\beq\label{defE}
\E(V)=\E(v_{1},v_{2})=\frac{1}{2}\left\|\pa_{x} V\right\|_2^2 
-\frac1{p+1}\int\big(|v_1|^{2p+2}+|v_2|^{2p+2}
+2\beta |v_1v_2|^{p+1}\big),
\eeq
when  the exponent $p$ satisfies
\beq\label{pzero}
1\leq p<2.
\eeq
The interest in finding ground states is  also motivated by their properties with respect of the analysis of the dynamical system \eqref{schr}, such as stability properties. For the single Schr\"odinger equation many notions of stability have been introduced and proved, among all, we recall \cite{cl} and \cite{weinsteinMS,weinsteinCpam}; in the former it is proved that the ground state, which is unique, of the equation 
\begin{equation}	\label{eeqr}
-\frac12\pa_{xx} z+z=z^{2p+1}\quad \text{in $\R$},
\end{equation}
is orbitally stable, that is, roughly speaking,  if $\phi^{0}$ is a function close to  $z$ with respect to the $H^{1}$ norm then the solution of the Cauchy problem
\begin{equation}\label{scalSE}
\begin{cases}
\dys \iu\pa_{t}\phi+ \frac12\pa_{xx}\phi+|\phi|^{2p}\phi=0 & \text{in $\R\times\R^{+}\!\!,$} \\
\phi(0,x)=\phi^0(x) & \text{in $\R,$}
\end{cases}
\end{equation}
where $\phi:[0,\infty)\times\R\to\C$, $\phi^{0}:\R\to\C$ and 
$1\leq p<2$, remains close to $z$ up to phase rotations and translations.
In \cite{weinsteinMS,weinsteinCpam} the study becomes deeper assuming that $z$ is non-degenerate, that is the linearized operator for~\eqref{eeqr}
has a $1$-dimensional kernel which is spanned by $\partial_{x}z$. More precisely, it is proved that for every $\phi\in H^1(\R)$ such that 
$\|\phi\|_{L^2}=\|z\|_{L^2}$, the following inequality holds
\begin{equation}\label{modscal}
\E(\phi)-\E(z)\geq C\inf_{x_{0}\in \R \atop\theta\in [0,2\pi)}
\|\phi-e^{i\theta}z(\cdot-x_0)\|_{H^1}^{2},
\end{equation}
for some positive constant $C$, provided that the energy $\E(\phi)$
is sufficiently close to $\E(z)$. Here, $\E$ is the energy defined in 
\eqref{defE} once we consider $V=(z,0)$. Inequality \eqref{modscal} allows to provide not only the same orbital stability result proved in \cite{cl}, but it also permits to derive explicit differential equation to which the phase and position adjustment have to obey for the ground state to be linearly stable. Moreover,  \eqref{modscal}
tells us that the energy functional can be seen as a Lyapunov functional, as it
measures the deviation of the solution of \eqref{schr} from the ground state orbit.
\\
The main goal of this paper is to extend inequality~\eqref{modscal} 
to the more general framework of 1D vector Schr\"odinger problems. 
In order to do this we are lead to consider non-degenerate ground state for system \eqref{ellittico}. This notion is introduced in the following definition.
\begin{definition}\label{nondeg}
We will say that a ground state solution $R=(r_1,r_2)$ of 
system~\eqref{ellittico} is non-degenerate if the set of 
solutions of the linearized system
\begin{equation}\label{linsyste}
\begin{cases}
-\frac{1}{2}\pa_{xx}\phi+\phi=[(2p+1)r_1^{2p}+
\beta p r_1^{p-1}r_2^{p+1}]\phi
+\beta(p+1)r_1^pr_2^p\psi & \text{in $\R$,} \\
-\frac{1}{2}\pa_{xx}\psi+\psi=[(2p+1)r_2^{2p}
+\beta p r_1^{p+1}r_2^{p-1}]\psi+
\beta(p+1)r_1^pr_2^p\phi & \text{in $\R$,}	
\end{cases}
\end{equation}
is an $1$-dimensional vector space and any solution $(\phi,\psi)$ 
of~\eqref{linsyste} is given by $\theta\pa_{x} R$, 
for some $\theta\in\R$.
\end{definition}

The main result of the paper is stated in the following

\bte\label{main}
Let $R$ be non-degenerate and assume~\eqref{pzero}. Then, for every $\Phi\in H^{1}\times H^{1}$ with 
$$
\|\Phi\|_{L^{2}\times L^{2}}=\|R\|_{L^{2}\times L^{2}},
$$ 
the following inequality holds
\begin{align*}
\E(\Phi)-\E(R)\geq &
\inf_{x\in \R \atop\theta\in [0,2\pi)^2}\|\Phi-(e^{\iu\theta_{1}}
r_{1}(\cdot-x),e^{\iu\theta_{2}}r_{2}(\cdot-x))\|_{H^{1}\times H^{1}}^{2}
\\
&\quad
+\text{{\small o}}\Big(
\inf_{x\in \R\atop\theta\in [0,2\pi)^2}\|\Phi-(e^{\iu\theta_{1}}
r_{1}(\cdot-x),e^{\iu\theta_{2}}r_{2}(\cdot-x))\|_{H^{1}\times H^{1}}^{2}\Big)
\end{align*}
where $\text{{\small o}}(x)$ satisfies 
$\text{ {\small o}}(x)/x\to 0$ as $x\to 0$.
\ete

As interesting consequences, we will obtain the property of being isolated, and of being orbitally stable for a non-degenrate ground state.
In \cite{mmp3} it has been recently proved that the set of ground states of \eqref{ellittico} enjoys the orbital stability property. To this respect, we have to  recall that up to now it is not yet been proved a uniqueness result for ground state solutions of the system \eqref{ellittico}. Therefore,  a solution of \eqref{schr} which starts near a ground state $R$, may leave the orbit around $R$ and approach the orbit generated by another ground state. But, this is not the case, once we know that the ground states are isolated. 
This property is easily obtained as a consequence of Theorem \ref{main} as stated in the following corollary.
\bco\label{cor1}
Let $R$ be non-degenerate and assume~\eqref{pzero}.
Then $R$ is isolated, that is, if there exists a  ground state 
of~\eqref{ellittico} $S$  satisfying $\|R-S\|_{\H}<\delta$ for a 
$\delta>0$ sufficiently small, then $S=R$ up to a translation 
and a phase change.
\eco
Then,  we can also prove the following
\bco
\label{cor2}
Let $R$ be non-degenerate and assume~\eqref{pzero}. Then $R$ is orbitally stable.
\eco
We recall that a ground  state $R=(r_{1},r_{2})$ is  
said to be orbitally stable if for  any given $\eps>0$,  
there exist $\dt(\eps)>0$ such that
$$
\sup_{t\in [0,\infty)}\inf_{x\in\R\atop \theta\in [0,2\pi)^2} 
\| \Psi(t,\cdot)-(e^{\iu\theta_{1}} r_{1}(\cdot-x),
e^{\iu\theta_{2}}r_{2}(\cdot-x)\|_{H^{1}\times H^{1}}< \eps,
$$
provided that
$$
\inf_{x\in\R\atop \theta\in [0,2\pi)^2} \| \Psi^{0}-(e^{\iu\theta_{1}}
r_{1}(\cdot-x),e^{\iu\theta_{2}}r_{2}(\cdot-x)\|_{H^{1}\times H^{1}}< \dt,
$$
where $\Psi$ is the solution of~\eqref{schr} with initial
datum $\Psi^{0}$.
\\
Theorem \ref{main} plays   a very important 
role also in the study of the so-called {\em soliton dynamics} for 
Schr\"o\-din\-ger. More precisely,  when one considers \eqref{schr}
when the Plank's constant $\hbar$ explicitly appears in the equations,
and studies  the evolution, in the semi-classical limit ($\hbar\to 0$),  of  the solution of \eqref{schr} starting from a $\hbar$-scaling of a soliton, once 
the action of external forces  appears.  We refer the reader to
\cite{bj,K1,K2} for the scalar case and to~\cite{mps-soliton} 
for systems, where the authors have recently showed, in
semi-classical regime, how the  soliton dynamics can be derived from 
Theorem \ref{main}. 
\\
Finally, we have to point out that some of our results can be proved in general dimension $n\geq 1$ as well, with minor changes. Unfortunately, this is not the case for our main Theorem, since, in order to work on the linearized equation, and to perform Taylor expansion on the energy functional ${\mathcal E}$,
we need enough regularity on the nonlinear term and this forces us to restrict the range of $p$ because of the presence of the coupling term.  
Of course, it is a really interesting open problem,  to prove the assertion of 
Theorem~\ref{main}   for any $n\geq 1$ and any $0<p<2/n$.
\vskip8pt
In Section~\ref{spectral}, we will study some delicate spectral 
properties of the linearized system introduced in Definition~\ref{nondeg}. 
The proofs of Theorem~\ref{main} and of Corollaries~\ref{cor1} and~\ref{cor2} will be carried out in Section~\ref{proofsection}.
Finally, in Section~\ref{ground}, we shall prove that there exists a
non-degenerate ground state for system~\eqref{ellittico}.

\section{Spectral analysis of the linearized operators}\label{spectral}

In this section we will prove some important properties concerning
the linearized Schr\"odinger system associated with~\eqref{schr}. 
\\
We will make use of the  functional spaces
$\L=\ldue(\R,\C)\times\ldue(\R, \C)$ and
$\H=\huno(\R,\C)\times\huno(\R,\C)$. 
We recall that the inner product between $u,\,v\in \C$ is given by 
$u\cdot v=\Re(u\bar v)=1/2(u\bar v+v\bar u)$. It is known (see \cite{caz,susu}) that \eqref{schr} is 
well locally posed in time, for any $p$, in the space $\H$ endowed with the norm 
$\|\Phi\|_{\H}^{2}=\|\pa_{x}\Phi\|_{2}^{2}+\|\Phi\|_{2}^{2}$ 
for every $\Phi=(\phi_1,\phi_2)\in\H$. Moreover we set
the ${\varmathbb{L}^{q}}$ norm as
$\|\Phi\|_q^q=\|\phi_1\|_q^q+\|\phi_2\|_q^q$ 
for any $q\in\left[1,\infty\right)$, we denote by
$(U,V)$ the inner scalar product in $\L$ and by $(U,V)_{\H}$ 
the inner scalar product in $\H$.
In \cite{fm} it is proved that, for $p$ satisfying $0<p<2$
the solution of the Cauchy problem \eqref{schr} exists globally 
in time and  the mass  of a solution and its total energy are 
preserved in time, that  is having defined the total energy of 
system~\eqref{schr} as
\beq\label{energy}
\E\left(\Phi(t)\right)=\frac{1}{2}\left\|\pa_{x}\Phi (t)\right\|_2^2 
-\int F\left(\Phi(t)\right) 
\eeq
where
\beq\label{defF}
F(U)=F(u_{1},u_{2})=\frac1{p+1}\Big(|u_{1}|^{2p+2}+|u_{2}|^{2p+2}
+2\beta |u_{1}u_{2}|^{p+1}\Big),
\eeq
the following conservation laws hold  (see \cite{fm}):
\beq\label{mass}
\|\phi_1\|_2^2=\|\phi^{0}_{1}\|_2^2,\qquad
\|\phi_2\|_2^2=\|\phi^{0}_{2}\|_2^2,
\qquad \E\left(\Phi(t)\right)=\E(0)=\frac{1}{2}\left\|\pa_{x}\Phi^{0}\right\|_2^2
-\int F\left(\Phi^{0}\right).
\eeq
Setting $\phi_i=r_i+\eps w_i$, $i=1,2$, the linearized 
Schr\"odinger system at $r_i$ in $w_i$ is given by
\beq\label{schrLIN}
\begin{cases}
\dys i\pa_{t}w_{1}+ \frac12\pa_{xx} w_{1}-w_1+G_1(w_1,w_2)=0 & \text{in $\R$}, 
\smallskip
\\ 
\dys i\pa_{t}w_{2}+\frac12\pa_{xx} w_{2}-w_2+G_2(w_1,w_2)=0 & \text{in $\R$},  
\end{cases}
\eeq
where we have set
$$
G_1(w_1,w_2)  =
\left[r_1^{2p}+\beta r_1^{p-1}r_2^{p+1}\right]w_1
+\left[2pr_1^{2p}+\beta(p-1)r_1^{p-1}r_2^{p+1}\right]\Re(w_1)
+\beta (p+1)r_1^pr_2^p \Re(w_2) ,
$$
$$
G_2(w_1,w_2) =
\left[r_2^{2p}+ \beta r_1^{p+1}r_2^{p-1}\right]w_2
+\left[2pr_2^{2p}+\beta(p-1)r_1^{p+1}r_2^{p-1}\right]\Re(w_2)
+\beta(p+1)r_1^pr_2^p\Re(w_1).                         
$$
System~\eqref{schrLIN} can be written down as
$\partial_t W=LW,$ for $L:\L\times \L\to \L\times \L$ defined by
\begin{equation*}
L=\left(
\begin{array}{cc}
0 & L_{-}
\\\\
-L_{+} & 0
\end{array}\right),\qquad W\in\C^2,\, W=(w_{1},w_{2})
\end{equation*}
and where the operators $L_{-},\,L_{+}:L^{2}(\R,\R)\times L^{2}(\R,\R)
\to L^{2}(\R,\R)\times L^{2}(\R,\R)$ acting respectively 
on the real and imaginary parts  of $w_{i}$.
 are the following
\beq
\label{operator}
L_{+}=\left(
\begin{array}{cc}
L_{+}^{11} & L_{+}^{12}
\\\\
L_{+}^{21} & L_{+}^{22}
\end{array}
\right)
\qquad L_{-}=\left(
\begin{array}{cc}
L_{-}^{11} & 0
\\\\
0 & L_{-}^{22}
\end{array}
\right)
\eeq
where $L_{+,-}^{ij}:L^{2}(\R,\R)\to L^{2}(\R,\R)$ are defined by
\begin{align*}
L_{+}^{11}&=\,-\frac12\pa_{xx} +1-H^{11}(R)
\qquad 
L_{+}^{12}=L_{+}^{21}=-H^{12}(R)
\\\noalign{\vskip4pt}
L_{+}^{22}&=\,-\frac12\pa_{xx} +1-H^{22}(R)
\\\noalign{\vskip4pt}
L_{-}^{11}&=\,-\frac12\pa_{xx} +1-\left[r_{1}^{2p}+\bt r_{1}^{p-1}r_{2}^{p+1}\right]
\qquad\qquad
L_{-}^{22}=\,-\frac12\pa_{xx} +1-\left[r_{2}^{2p}+\bt r_{1}^{p+1}r_{2}^{p-1}\right]
\end{align*}
and the Hessian matrix $H_{F}(U)=(H^{ij}):(\R^{+})^{2}\to {M_{2\times 2}}(\R)$ is given by
\begin{align*}
H^{11}&= (2p+1)u_{1}^{2p}+p\bt u_{1}^{p-1}u_{2}^{p+1}
\qquad H^{12}=H^{21}=\, (p+1)\bt u_{2}^{p}u_{1}^{p}
\\\noalign{\vskip4pt}
H^{22}&= (2p+1)u_{2}^{2p}+p\bt u_{2}^{p-1}u_{1}^{p+1}.
\end{align*}

We will study $L_{+}$  on $ \mathcal{V}$, namely the closed subspace of $\H$ defined as
\beq\label{defV}
\mathcal{V}=\left\{U\in \H\,:\,(U,R)=0\right\}.
\eeq
The first important property of $L_{+}$ on $\mathcal{V}$ is proved in the following proposition. 
\bpr\label{L+prima}
Assume \eqref{pzero} and that $R$ a ground state of \eqref{ellittico}. Then
$\dys \inf\limits_{\mathcal{V}}\left(L_{+}(U),U\right)=0$. 
\epr
\bdim
First notice that $U_{*}=(r_1',r_2')$
belongs to $\mathcal{ V}$ and $U_{*}$
satisfies $(L_{+}(U_{*}),U_{*})=0$, showing that
the infimum is less or equal than zero.
On the other hand, since $R$ solves problem \eqref{elle2}, 
of course $R$ is also a minimum point of  $\I=\E(\Phi)+ \|\Phi\|_{2}^{2}$ 
on $\M$. Consequently, for any smooth curve $\varphi:[-1,1]\to \M$ such that $\varphi(0)=R$, it follows 
$$
\frac{d^{2}\I(\varphi(s))}{ds^{2}}\Bigg|_{s=0}\geq 0.
$$
Therefore, taking into account that $\I'(R)=0$, we get
\begin{align*}
0 &
 \leq \langle \I{''}(\varphi(s))\varphi'(s), \varphi'(s)\rangle\Big|_{s=0} 
 +\langle \I'(\varphi(s)),\varphi''(s)\rangle  \Big|_{s=0} 
\\& =\langle \I{''}(R)\varphi'(0),\varphi'(0)\rangle +\langle \I'(R),\varphi''(0)\rangle
=\langle \I{''}(R)\varphi'(0),\varphi'(0)\rangle.
\end{align*}
Now, taking into account that the map $s\mapsto \|\varphi(s)\|_{2}$ is constant, it readily follows that
$\varphi'(0)$ belongs to $\mathcal{ V}$, which yields the assertion by the arbitrariness of $\varphi$.
\edim

The above result is the first step to show that $L_{+}$ is coercive once we restrict it on a closed subspace of $\mathcal{V}$, as shown in the following proposition.
\bpr
\label{oneofmaitoools}
Assume \eqref{pzero} and that $R$ is a ground state of \eqref{ellittico}
satisfying Definition \ref{nondeg}. Then
\beq\label{infpos}
\inf_{U\in\mathcal{V}_{0},\;\;  }
\frac{\left(L_{+}(U),U\right)}{\|U\|^{2}_{2}}>0,\qquad
\mathcal{V}_{0}=\left\{U\in\H : (U,R)=(U,H_{F}(R)\pa_{x}R)=0\right\}.
\eeq

\epr
\bdim
Denoting with $\alpha$ the infimum 
$$
\alpha=\inf_{\|V\|_{L^2}=1,\,V\in\mathcal{V}_{0}}(L_{+}(V),V),
$$
first notice that Proposition \ref{L+prima} implies that $\alpha$ is 
nonnegative, so that we only have to show that $\alpha$ is not zero.
Let us argue by contradiction and suppose that $\alpha=0$. Taken $U_{n}$  a minimizing sequence, 
from the regularity properties of $R$ it follows that ${U}_{n}$ is bounded in $\H$. 
These gives us a function $U\in \H$, such that  $U_{n}\rightharpoonup U$ weakly (up to a subsequence) in $\H$, implying that  $U\in \mathcal{V}_{0}$. From Proposition \ref{L+prima} and
\eqref{infpos}, we get
$$
0\leq \left(L_{+}(U),U\right)  \leq \liminf_{n\to\infty} 
\left\{\|U_{n}\|^{2}_{\H}-(U_{n},H_{F}(R)U_{n})\right\} 
=\lim_{n\to\infty} \left(L_{+}(U_{n}),U_{n}\right)=0.
$$
So that $U$ solves $\left(L_{+}(U),U\right)=0$ and
$\left(L_{+}(U_{n}),U_{n}\right)\to \left(L_{+}(U),U\right)$. Moreover,
\begin{align*}
\|U\|^{2}_{\H} \leq \liminf_{n\to\infty}\|U_{n}\|^{2}_{\H}
& \leq \limsup_{n\to\infty}\|U_{n}\|^{2}_{\H}= \lim_{n\to\infty} 
\big\{\left(L_{+}(U_{n}),U_{n}\right) +(U_{n},H_{F}(R)U_{n})\big\}\\
& =\left(L_{+}(U),U\right)+(U,H_{F}(R)U)=\|U\|^{2}_{\H},
\end{align*}
from which $U_{n}\to U$ strongly in $\H$, so that $\|U\|_{\L}=1$ and $U$ solves  the constrained minimization problem \eqref{infpos}. 
When we derive the functional $(L_{+}(V),V)/\|V\|_{\L}^{2}$ and use
that $(L_{+}(U),U)=0$ we obtain that there exists  Lagrange multipliers
$\mu,\,\gamma\in \R$ such that 
\beq\label{lagrange} 
\left(L_{+}U,V\right)=\mu\left( R,V\right)+\left(\gamma\cdot H_{F}(R) \pa_{x} R,V\right),\qquad \text{for every $V\in \H$}.
\eeq
Choosing as test function $V=\pa_{x} R$ and taking into consideration that $(R,\partial_jR)=0$, gives
$$
0=\left(L_{+}(U),\pa_{x} R\right)=
\left(\gamma\cdot H_{F}(R)\pa_{x} R,\pa_{x} R\right)
=\gamma (H_{F}(R)\pa_{x}R,\pa_{x}R),
$$
where we have taken into account that $L_{+}$ is a self-adjoint 
operator and $\pa_{x}R=\left(\pa_{x}r_{1},\pa_{x}r_{2}\right)$ 
is a solution of $L_{+}V=0$. Since $R$ has even components the 
summands on the right hand side are nonzero, so that $\gamma=0$. 
As a consequence, $U$ solves $L_{+}U=\mu R$.  Moreover,
we consider the vector $x \cdot\pa_{x} R$, whose components are
$x \cdot\pa_{x} R=(x\pa_{x}r_{1},x\pa_{x}r_{2})$ and we compute
$L_{+}(x \cdot\pa_{x} R)$. After some simple calculations, one reaches
$$
L_{+}(x \cdot\pa_{x} R)=(-\pa_{xx} r_{1},-\pa_{xx} r_{2})
\quad
\text{and}
\quad
L_{+}(R/p)=-2(r_1^{2p+1}+\beta r_2^{p+1}r_1^{p},
r_2^{2p+1}+\beta r_1^{p+1} r_2^{p}).
$$
Then, in turn, we get $L_{+}(R/p+x \cdot\pa_{x} R)=-2R$, and by linearity 
$$
L_{+}\left(-\mu/2(R/p+x \cdot\pa_{x} R)\right)=\mu R.
$$
Then, Definition~\ref{nondeg} (nondegeneracy) immediately yields
\begin{equation}\label{Uequatt}
U=-\mu/2(R/p+x \cdot\pa_{x} R)+\theta\cdot\pa_{x} R
\end{equation} 
for some constant $\theta\in\R$. Now we have to show that 
$\theta=0$, by using the available constraints. By applying 
to equation~\eqref{Uequatt} the self-adjoint operator $H_F=H_F(R)$, we get
$$
H_FU=-\frac{\mu}{2p} H_F R-\frac{\mu}{2}H_F x\cdot\pa_{x} R+H_F\theta\cdot\pa_{x} R.
$$
As $U\in\mathcal{V}_{0}$, it results $
(H_FU,\pa_{x} R)=(U,H_F\pa_{x} R)=0.$
Furthermore, since $R$ is a radial solution of \eqref{ellittico}, 
we also have that $(H_F R,\pa_{x} R)= 
(H_F x\cdot\pa_{x} R,\pa_{x} R)=0.$ On the other hand
$$
(H_F\theta\cdot\pa_{x} R,\pa_{x} R)=\theta (H_F\pa_{x} R,
\pa_{x} R)=c\theta
$$
with $c\neq 0$, so it has to be $\theta=0$. 
Then~\eqref{Uequatt} reduces to
$$
U=-\frac{\mu}{2p}R-\frac{\mu}{2}x\cdot\pa_{x} R.
$$
Computing the $L^{2}$-scalar product with $R$ and keeping in 
mind that $U\in  \mathcal{V}_{0}$ yields
$$
0=(U,R)=-\frac{\mu}2\left[\frac1{p}\|R\|_{2}^{2}+(x\cdot\pa_{x} R,R)\right].
$$
As far as concern the last term in the previous relation, we integrate by parts
and obtain
$$
(x\cdot\pa_{x} R,R)=-\frac{1}2 \|R\|_{2}^{2}.
$$
The last two equations and \eqref{pzero} give the desired contradiction.
\edim

\bos\label{coerc}
The argument in the proof of the previous Proposition shows that
there exists a positive constant $\alpha_0$ such that
\beq
\label{coercL+}
(L_{+}V,V)\geq \alpha_{0} \|V\|_{2}^{2},\qquad\text{for all $V\in \mathcal{V}_{0}$}.
\eeq
Moreover, if we consider $|||U|||=\sqrt{(L_{+}U,U)}$  
for every $U\in \mathcal{V}_{0}$, we obtain that $|||\cdot|||$ 
satisfies all the required properties of a norm, by \eqref{coercL+} 
and by the self-adjointness property of $L_{+}$. 
In addition, every Cauchy sequence  $\{U_{n}\}$ with respect to
$|||\cdot|||$ has a  strong limit $U$ belonging $L^{2}$; moreover 
$U$ satisfies all the orthogonality relations required in 
$\mathcal{V}_{0}$. Besides, computing $(L_{+}(U_{n}-U_{m}),U_{n}-U_{m})$  
gives that also $\{\pa_{x} U_{n}\}$ is a Cauchy sequence in $L^{2}$ 
then $U$ is necessarily the strong limit of $\{U_{n}\}$ in $\H$. 
Finally, $|||U_{n}-U|||\to 0$ by the definition of $L_{+}$. 
As a consequence, $\mathcal{V}_{0}$ is a Banach space with respect 
to this norm, and we get the equivalence with the standard $\H$ 
norm, namely there exists $\alpha>0$ such that
$$
\left(L_{+}V,V\right)\geq \alpha \|V\|_{\H}^{2},\qquad\text{for all $V\in \mathcal{V}_{0}$}.
$$
\eos
Before stating our next result let us prove the following lemma.
\ble
Let us take $\Phi\in \L$ such that $\|\Phi\|_{2}=\|R\|_{2}$ and
consider the difference $W=\Phi-R$. Denoting with
$U$ and $V$ the real and imaginary part of $W$, it results
\beq\label{ideU}
\left(R,U\right)=-\frac12 \left[\|U\|_{2}^{2}+\|V\|_{2}^{2}\right]=-\frac12
\|W\|_{2}^{2}
\eeq
\ele
\bdim
The above identity immediately follows by imposing 
$\|R+W\|_{2}^{2}=\|R\|_{2}^{2}$ and by recalling that $R$ 
is a real function.
\edim

\bpr\label{fineL+}
Assume \eqref{pzero} and that $R$ satisfies Definition \ref{nondeg}. 
Moreover, let us take $W=U+iV$ satisfying \eqref{ideU} with $U$ verifying
\beq\label{ortder}
(U,H_{F}(R)\pa_{x}R)=0.
\eeq
Then, there exists positive constants $D,\,D_{i}$ such that
\beq\label{disL+}
\left(L_{+}U,U\right)\geq D\|U\|_{\H}^{2}-D_{1} \|W\|_{2}^{4}-D_{2}
\|W\|_{2}^{2}\|\pa_{x} W\|_{2}
\eeq
\epr
\bdim
Without loss of generality, we can suppose that $\|R\|_{2}=1$; moreover,  we decompose $U$ as  $U=U_{||}+U_{\bot}$
where $U_{||}=\left(U,R\right)R$, while $U_{\bot}=U-U_{||}$ is orthogonal to $R$ with respect to the $L^{2}$ scalar product. Since $L_{+}$ is self-adjoint it results
\beq\label{deco}
\left(L_{+}U,U\right)=\left(L_{+}U_{||},U_{||}\right)+2\left(L_{+}U_{\bot},U_{||}\right)
+\left(L_{+}U_{\bot},U_{\bot}\right).
\eeq
Next, we study separately the summands on the right hand side of this 
formula. Observe that, taking into account identity~\eqref{ideU}, we have 
\begin{equation}
	\label{graddineqq}
\|\pa_{x} U_\bot\|_{2}^2\geq \|\pa_{x} U\|_{2}^2-C\|W\|_{2}^{2}\|\pa_{x} W\|_{2},
\end{equation}
for some positive constant $C$. Since $(U_{||},H_{F}(R)\pa_{x}R)=0$, 
condition~\eqref{ortder} implies that also $U_{\bot}$ has to be orthogonal 
to $H_{F}(R)\pa_{x}R$, hence $U_{\bot}$ is in $\mathcal{V}_{0}$.
Then Remark~\ref{coerc},~\eqref{graddineqq} and~\eqref{ideU} give us 
\begin{align}\label{uort}
\left(L_{+}U_{\bot},U_{\bot}\right)& \geq D \|U_{\bot}\|_{\H}^{2}\geq
D\| U\|_{\H}^2-CD\|W\|_{2}^{2}\|\pa_{x} W\|_{2}-
D\|U_{||}\|_{2}^{2}\\
&=D\|U\|_{\H}^{2}-d_{1}\|W\|_{2}^{2}\left[ 
\|W\|_{2}^{2}+\|\pa_{x} W\|_{2}\right] . \notag
\end{align}
We also obtain from \eqref{ideU} that
\beq\label{prod}
\left(L_{+}U_{\bot},U_{||}\right)=(R,U)\left(L_{+}U_{\bot},R\right)
= -\frac12 \|W\|_{2}^{2}
\left(L_{+}U_{\bot},R\right)\geq
-d_{2} \|W\|_{2}^{2}\|\pa_{x} W\|_{2}.
\eeq
As far as concern the last term in \eqref{deco}, it results
$$
\left(L_{+}U_{||},U_{||}\right)=(U,R)^{2}\left(L_{+}R,R\right)
=\frac14\|W\|_{2}^{4}\left(L_{+}R,R\right)\geq -d_{3}\|W\|_{2}^{4}.
$$
This last equation, joint with \eqref{uort} and \eqref{prod} yields the conclusion.
\edim
\bpr
It results 
$\dys \inf\limits_{V\neq 0, \; (v_{i},r_{i})_{H^{1}}=0}
\frac{\left(L_{-}(V),V\right)}{\|V\|^{2}_{2}}>0$. 
\epr
\bdim
Let us first prove that $L_{-}$ is a positive operator.
Denoting with $\sigma_{d}(L_{-})$ the discrete spectrum of the operator $L_{-}$ it results 
\beq\label{spettro}
\sigma_{d}(L_{-})=\sigma_{d}(L^{11}_{-})\cup\sigma_{d}(L^{22}_{-}).
\eeq 
Indeed, if $\lambda\in \sigma_{d}(L^{11}_{-})$ we get that
$L^{11}_{-}(u)=\lambda u$, then $\lambda\in \sigma_{d}(L_{-})$ with eigenfunction $U=(u,0)$, analogous argument holds for $\lambda\in \sigma_{d}(L^{22}_{-})$, proving that $\sigma_{d}(L^{11}_{-})\cup\sigma_{d}(L^{22}_{-})\subseteq \sigma_{d}(L_{-})$.
On the other hand, if $\lambda\in \sigma_{d}(L_{-})$ there exists
$U=(u_{1},u_{2})\neq (0,0)$ such that 
$$
L_{-}^{11}u_{1}=\lambda u_{1},\quad 
L_{-}^{22}u_{2}=\lambda u_{2}
$$
so that, if $u_{1}\neq 0$ $\lambda\in \sigma_{d}(L^{11}_{-})$, otherwise
$u_{2}\neq 0$ and  $\lambda\in \sigma_{d}(L^{22}_{-})$, showing 
\eqref{spettro}.
Moreover, since $L_{-}R=0$, with $R=(r_{1},r_{2})\neq (0,0)$, 
$r_{i}\geq 0$, we get that $\lambda =0$ is the first eigenvalue of 
$L^{11}_{-}$ and $L^{22}_{-}$ when both $r_{1},\,r_{2}\neq 0$.
Besides, if for example $r_{1}\equiv 0$, $\lambda=0$ is the first 
eigenvalue of $L_{-}^{22}$, while $L_{-}^{11}=-\pa_{xx}+1$ and its 
discrete spectrum  is empty (see e.g. Chapter 3 in \cite{bershu}),  
yielding that $\lambda=0$ is the first eigenvalue of $L_{-}$.
Then $(L_{-}(V),V)\geq 0$ for every function $V\in\H$, proving that 
$L_{-}$ is a positive operator. Arguing now as in the proof of
Proposition~\ref{oneofmaitoools}, and considering the (nonnegative) infimum
$$
\alpha=\inf_{\|V\|_{L^2}=1,\,\, (V_{i},r_{i})_{H^{1}}=0}(L_{-}(V),V),
$$
assuming by contradiction that $\alpha=0$, we find that there exists a nonzero 
minimizer $U$ (satisfying the constraints) for the problem such that
\begin{equation}
	\label{minimizconstrain}
(L_{-}U,U)=0
\end{equation}
Taking into account that the constraints $(U_{i},r_{i})_{H^{1}}=0$ can be written in the $L^2$ form
\begin{equation}
	\label{constrains-i}
(q^{11}_{-}(R)R_1,U)=0,\qquad (q^{22}_{-}(R)R_2,U)=0,
\end{equation}
where we have set
$$
q^{11}_{-}(R)=r_{1}^{2p}+\bt r_{1}^{p-1}r_{2}^{p+1},
\qquad
q^{22}_{-}(R)=r_{2}^{2p}+\bt r_{1}^{p+1}r_{2}^{p-1},
\quad
R_1=(r_1,0),
\quad
R_2=(0,r_2).
$$
we have three lagrange parameters $\lambda,\gamma_1,\gamma_2\in\R$ such that
$$
(L_{-}U,V)=\lambda (U,V)+\gamma_1(q^{11}_{-}(R)R_1,V)+\gamma_2(q^{22}_{-}(R)R_2,V)
$$
for all $V\in\H$. Hence, by choosing $V=U$ and taking into account~\eqref{minimizconstrain}
and that $U$ satisfies the constraints~\eqref{constrains-i}, we immediately get $\lambda=0$.
Choosing now $V=R_1$ and $V=R_2$ and taking into account $L_{-}$ is self-adjoint and
that $L_{-}R_i=0$ we obtain $\gamma_1=\gamma_2=0$. Therefore, we conclude that
$$
L_{-}U=0,
$$
namely $L_{-}^{11}u_1=0$ and $L_{-}^{22}u_2=0$ where we set $U=(u_1,u_2)$. In turn,
$u_i$ is a first eigenfunction of $L_{-}^{ii}$, which yields $u_i\in{\rm span}(r_i)$
since the first eigenvalue is simple (see e.g.\ Theorem 3.4 in~\cite{bershu}). 
This is of course a contradiction with~\eqref{constrains-i}.
Hence $\alpha>0$ and the proof is complete.
\edim
\bos
\label{stimaLmeno}
Arguing as in Remark~\ref{coerc}, it is possible to find 
a positive constant $\alpha>0$ such that
\begin{equation*}
(L_{-}V,V)\geq \alpha \|V\|_{\H}^2,\qquad 
\text{for all $V\in\H$ with $(v_{i},r_{i})_{H^{1}}=0$, $i=1,2$}.
\end{equation*}
\eos


\section{Proofs of the main results}
\label{proofsection}

In order to prove Theorem \ref{main}, the following characterization will be crucial.
\bpr
\label{ortogonal}
Let us consider $y_{0}\in \R$ and $\Gamma=(\gamma_1,\gamma_2)\in \rdue$ be such that 
\beq\label{minimo}
\min_{x_{0}\in \R \atop\Theta\in \rdue}
\|(\phi_{1}(\cdot+x_{0})e^{\iu\theta_{1}},\phi_{2}(\cdot+x_{0})e^{\iu\theta_{2}})-R\|_{\H}^{2}
=
\|(\phi_{1}(\cdot+y_{0},t)e^{\iu\gamma_{1}},\phi_{2}(\cdot+y_{0})
e^{\iu\gamma_{2}})-R\|_{\H}^{2}
\eeq
Then, writing 
$$
(\phi_{1}(\cdot+y_{0},t)e^{\iu\gamma_{1}},\phi_{2}(\cdot+y_{0},t)
e^{\iu\gamma_{2}})=R+W,
$$ 
where $W=U+\iu V$, the following orthogonality condition are satisfied
\beq\label{orto}
\left(U, H_{F}(R)\pa_{x} R \right)=0, \qquad
\left(v_{1},r_{1}\right)_{H^{1}}=\left(v_{2},r_{2}\right)_{H^{1}}=0.
\eeq
\epr
\bdim
Let us introduce the functions $P,\,Q:\R\times \rdue\to \R$ 
defined by
\begin{align*}
P(x_{0}, \Theta) &=P(x_{0}, \theta_{1},\theta_{2})=\|(\phi_{1}(\cdot+x_{0})e^{\iu\theta_{1}},
\phi_{2}(\cdot+x_{0})e^{i\theta_{2}})-R\|_{2}^{2}
\\
Q(x_{0}, \Theta) &=Q(x_{0}, \theta_{1},\theta_{2})=\|
(\pa_{x} \phi_{1}(\cdot+x_{0})e^{\iu\theta_{1}},\pa_{x} 
\phi_{2}(\cdot+x_{0})e^{\iu\theta_{2}})-\pa_{x} R\|_{2}^{2}.
\end{align*}
Writing down the partial derivatives of $P$ and $Q$ and 
integrating by parts, give us
\begin{align*}
\pa_{x_{0}}P(x_{0},\Theta) & =\sum_{j=1}^{2}\int\left(
\phi_{j}e^{\iu\theta_{j}}-r_{j}\right)e^{-\iu\theta_{j}}
\pa_{x_{0}}\overline{\phi}_{j}+
\left(\overline{\phi}_{j}e^{-\iu\theta_{j}}-r_{j}\right)e^{\iu\theta_{j}}
\pa_{x_{0}}\phi_{j}
\\
&=-2\sum_{j=1}^{2} \int r_{j}\Re\left(e^{\iu\theta_{j}}
\pa_{x_{0}} \phi_{j}\right);
\\
\pa_{x_{0}} Q(x_{0}, \Theta)  & =\sum_{j=1}^{2}\int
\pa_{x}\left(\phi_{j}e^{\iu\theta_{j}}-r_{j}\right) 
\pa_{x} \pa_{x_{0}} \overline{\phi}_{j} e^{-\iu\theta_{j}}+
\pa_{x}\left(\overline{\phi}_{j}e^{-\iu\theta_{j}}-r_{j}\right) 
\pa_{x} \pa_{x_{0}} \phi_{j} e^{\iu\theta_{j}}
\\
&=- 2\sum_{j=1}^{2} \int \pa_{x} r_{j}\Re\left(\pa_{x} 
\pa_{x_{0}}\phi_{j} e^{\iu\theta_{j}}\right);
\\
\frac{\partial P}{\partial \theta_{j}}(x_{0}, \Theta)  &=\iu\int\left[
-\left(\phi_{j}e^{\iu\theta_{j}}-r_{j}\right) e^{-\iu\theta_{j}} \overline{\phi}_{j}+
\left(\overline{\phi}_{j}e^{-\iu\theta_{j}}-r_{j}\right) e^{\iu\theta_{j}}  \phi_{j}\right]
\\
&=2 \int r_{j}\Im\left(e^{\iu\theta_{j}}\phi_{j}\right);
\\
\frac{\partial Q}{\partial\theta_{j}} (x_{0}, \Theta) &=\iu\int\left[
-\pa_{x}\left(\phi_{j}e^{\iu\theta_{j}}-r_{j}\right) \pa_{x}\overline\phi_{j} e^{-\iu\theta_{j}}+
\pa_{x}\left(\overline{\phi}_{j}e^{-\iu\theta_{j}}-r_{j}\right) \pa_{x}\phi_{j} e^{\iu\theta_{j}}\right]
\\
&=2\int \pa_{x} r_{j}\Im\left(\pa_{x} \phi_{j} e^{\iu\theta_{j}}\right).
\end{align*}
If $x_{0}=y_{0}$ and $\Gamma=(\gamma_1,\gamma_2)$ realize 
the minimum in~\eqref{minimo}, the following equations 
are satisfied
\begin{align*}
\frac{\partial (P+Q)}{\partial x_{0}}(x_{0}, \Theta) &=
-2\sum_{j=1}^{2} \int \left[r_{j}(x)\Re\left(e^{\iu\gamma_{j}}\frac{\partial {\phi}_{j}}{\partial x_{0}}(x-y_{0})\right)+\pa_{x} r_{j}(x) \Re\left(e^{\iu\gamma_{j}}\pa_{x} \frac{\partial {\phi}_{j}}{\partial x_{0}} (x-y_{0}) \right)\right]=0 
\\
\frac{\partial (P+Q)}{\partial \theta_{j}}(x_{0}, \Theta) &=
2\int\left[ r_{j}(x)\Im\left(e^{\iu\gamma_{j}}\phi_{j}(x-y_{0})\right)
+\pa_{x} r_{j}(x)\Im\left(e^{\iu\gamma_{j}}\pa_{x} \phi_{j}(x-y_{0}) \right)\right]=0.
\end{align*}
Denoting with $U$ and $V$ the real and imaginary (respectively) part
of $W=\Phi(x-y_{0})e^{\iu\Gamma}-R(x)$ and taking into account that
$R$ is real and does not depend on $x_{0}$, it follows
\begin{align*}
\frac{\partial (P+Q)}{\partial x_{0}}(x_{0}, \Theta) &=
\sum_{j=1}^{2} \int \left[r_{j}\frac{\partial u_{j}}{\partial x_{0}}+\pa_{x} r_{j}\pa_{x} \frac{\partial u_{j}}{\partial x_{0}}\right]
=-\sum_{j=1}^{2} \int \left[u_{j}\frac{\partial r_{j}}{\partial x_{0}}+\pa_{x} u_{j}\pa_{x} \frac{\partial r_{j}}{\partial x_{0}}\right]=0 
\\
\frac{\partial (P+Q)}{\partial \theta_{j}}(x_{0}, \Theta) &=
\int\left[ r_{j} v_{j}+\pa_{x} r_{j}\pa_{x} v_{j}\right]=0,\quad j=1,2.
\end{align*}
The second line of the above equations can be read as the orthogonality 
conditions on $V$ in~\eqref{orto}. As far as regards $U$, we only have 
to notice that $\pa_{x}R$ satisfies the linearized system 
of~\eqref{ellittico} so that all the conditions in~\eqref{orto} are proved.
\edim

\noindent
We are now ready to complete the proof of the main result, Theorem~\ref{main}.
\vskip4pt
\noindent{\bf Proof of Theorem~\ref{main} concluded.}
Let us consider $\Phi\in \H$ with $\|\Phi\|_{2}=\|R\|_{2}$ and 
$W(x)=\Phi(x-y_{0})e^{\iu\Gamma}-R(x)$, where $y_{0}\in\R$ and $\Gamma\in\R^2$
satisfy the minimality conditions~\eqref{minimo}. We want to control the $\H$ norm of $W$ in terms 
of the difference $\I(\Phi)-\I(R)$, being $\I$ is the action functional associated 
to the system and defined as
$$
\I(\Phi)=\E(\Phi)+ \|\Phi\|_{2}^{2}.
$$
To this aim, we first compute the difference $\I(\Phi)-\I(R)$ and we use scale invariance, obtaining
$\I(\Phi)-\I(R)=\I(R+W)-\I(R)$. Then, recalling that $\langle \I'(R),W\rangle=0$, Taylor expansion gives
\begin{align*}
\I(\Phi)-\I(R)&=\I(R+W)-\I(R)=\langle \I'(R),W\rangle+\langle \I''(R+\vth W)W,W\rangle  \\
&=\langle \I''(R)W,W\rangle +\langle \I''(R+\vth W)W,W\rangle-\langle \I''(R)W,W\rangle.
\end{align*}
In order to evaluate the difference on the right hand side we will use the $C^{2}$ regularity of $\I$, at this point it is crucial \eqref{pzero}. 
For simplicity, let  us consider separately the nonlinear terms in $\I$. 
The term $G:\H\to \R$ defined by
$$
G(U)=G(u_1,u_2)= \|u_1\|_{2p+2}^{2p+2}+\|u_2\|_{2p+2}^{2p+2},
$$
is of class $C^{3}$, as $p\geq 1$, so that 
\beq\label{gsec}
\langle G''(R+\vth W)W,W\rangle-\langle G''(R)W,W\rangle\geq -c_{1}\|W\|_{\H}^{3}.
\eeq
As far as concern the coupling term $\Upsilon:\H\to \R$ defined by
$\Upsilon(U)=\Upsilon(u_1,u_2)=\|u_1u_{2}\|_{p+1}^{p+1},$
it results
\begin{align*}
\langle \Upsilon''(U)W,W\rangle & =
(p^2-1)\int |u_{1}|^{p-3}|u_2|^{p-3}\left[
|u_{2}|^{4}\Re^2(u_1)|w_1|^2 +|u_{1}|^{4}\Re^2(u_2)|w_2|^2
\right]
\\
&
+(p+1)\int  |u_1|^{p-1}|u_2|^{p-1}
\left[|u_2|^{2}|w_1|^2 +|u_1|^{2}|w_2|^2 \right]
\\
&+2(p+1)^2\int |u_1|^{p-1}|u_2|^{p-1}\Re(u_1)\Re(u_2)\Re(w_1\overline{w}_2).
\end{align*}
When we write the difference $\langle \Upsilon''(R)W,W\rangle-
\langle \Upsilon''(R+\vth W)W,W\rangle $ we use that $R$ is a real 
function and we 
control the first two terms with the real parts by the modulus; finally we use 
the inequality
$$
\left||r_j+\vth w_j|^{p-1}-|r_j|^{p-1}\right|\leq C|w_j|^{p-1},
$$
to get 
\beq
\langle \Upsilon''(R)W,W\rangle-
\langle  \Upsilon''(R+\vth W)W,W\rangle
\geq -c_{1}\|W\|_{\H}^{2+\mu}\qquad\text{for some $\mu>0$.}
\eeq 
This inequality joint with \eqref{gsec} implies that
\beq
\langle \I''(R+\vth W)W,W\rangle
-\langle \I''(R)W,W\rangle
 \geq -C \|W\|_{\H}^{2+\mu}.
\eeq 
Therefore,
\begin{equation*}
\I(\Phi)-\I(R)\geq \langle \I''(R)W,W\rangle-C\|W\|_{\H}^{2+\mu}=
\langle L_- V,V\rangle+\langle L_+ U,U\rangle-C\|W\|_{\H}^{2+\mu}.
\end{equation*}
Taking into account the orthogonality conditions
of Proposition~\ref{ortogonal}, the assertion now follows from Proposition~\ref{fineL+}
and Remark~\ref{stimaLmeno}.
\edim

\noindent{\bf Proof of Corollary~\ref{cor1}}\quad
Let $\delta$ be a positive number to be chosen later. 
Moreover, let $R=(r_1,r_2)\in\H$ and $S=(s_1,s_1)\in \H$ be two given 
non-degenerate ground state solutions to system~\eqref{ellittico}  such that
$$
\|R-S\|_{\H}^2<\delta.
$$
Then, taking into account
the variational characterization~\eqref{elle2} for ground states, we learn that
$$
\E(R)=\E(S),\qquad \|R\|_{\L}=\|S\|_{\L}.
$$
Notice also that
$$
\inf_{x_{0}\in \R\atop\theta\in \rdue}\|R-(e^{i\theta_{1}}
s_{1}(\cdot-x_{0}),e^{i\theta_{2}}s_{2}(\cdot-x_{0}))\|_{\H}^{2}
\leq \|R-S\|_{\H}^2<\delta.
$$
Therefore, by applying Theorem~\ref{main}, if $\delta>0$ is chosen sufficiently small, we get
\begin{equation*}
\inf_{x_{0}\in \R \atop\theta\in \rdue}\|R-(e^{i\theta_{1}}
s_{1}(\cdot-x_0),e^{i\theta_{2}}s_{2}(\cdot-x_0))\|_{\H}^{2}\leq 0.
\end{equation*}
In turn we conclude that $R=S$, up to a suitable translation and phase change.
\vskip4pt
\noindent

\noindent{\bf Proof of Corollary~\ref{cor2}}\quad
Let $T>0$ and let us fix $\eps>0$ sufficiently small. 
Consider the solution $\Psi$ of system~\eqref{schr} with initial datum $\Psi^{0}$.
By the conservation laws, we have
$$
\|\Psi(t)\|_{\L}=\|\Psi^0\|_{\L},\,\,\quad   \E(\Psi(t))=\E(\Psi^0),\quad
\text{for all $t\in [0,\infty)$}.
$$
By the continuity of the energy $\E$, there exists $\dt=\dt(\eps)>0$ such that
$$
\E(\Psi(t))-\E(R)=\E(\Psi^0)-\E(R)<\eps,\quad\text{for all $t\in [0,\infty)$},
$$
provided that 
\begin{equation}
	\label{deltasmall}
\inf_{\theta\in\rdue
\atop x\in\R} \| \Psi^0(\cdot)-(e^{i\theta_{1}}
r_{1}(\cdot-x),e^{i\theta_{2}}r_{2}(\cdot-x))\|_{\H}^2< \delta.
\end{equation}
Then, if we define for any $t>0$ the positive number
\begin{equation*}
\Gamma_{\Psi(t)}=\inf_{\theta\in\rdue\atop x\in\R}
\|\Psi(t)-(e^{i\theta_1}r_1(\cdot-x),e^{i\theta_2}r_2(\cdot-x))\|_{\H}^{2},
\end{equation*}
we learn from Theorem~\ref{main} that there exist two positive constants ${\mathcal A}$ and $C$ such that
\begin{equation}
		\label{strongconclusion}
\Gamma_{\Psi(t)}\leq C(\E(\Psi(t))-\E(R)),
\end{equation}
provided that $\Gamma_{\Psi(t)}<{\mathcal A}$. Let us define the value
$$
T_0:=\sup\big\{t\in [0,T]:\,\, \Gamma_{\Psi(s)}<{\mathcal A}\,\,\text{for all $s\in [0,t)$}\big\}.
$$
Of course, it holds $T\geq T_0>0$ by means of~\eqref{deltasmall} (up to reducing the size of $\delta$,
if necessary) and the continuity of $\Psi(t)$.
Hence, we deduce that
\begin{equation}
	\label{firstconcl}
\sup_{t\in [0,T_0]}\inf_{\theta\in\rdue
\atop x\in\R} \| \Psi(t,\cdot)-(e^{i\theta_{1}}
r_{1}(\cdot-x),e^{i\theta_{2}}r_{2}(\cdot-x))\|_{\H}^2\leq 
C(\E(\Psi(t))-\E(R))=C(\E(\Psi^0)-\E(R))< C\eps .
\end{equation}
On the other hand, it is readily seen that, from this inequality, 
one obtains $T_0=T$. In fact, assume by contradiction
that $T_0<T$. Then, since by~\eqref{firstconcl}
$$
\Gamma_{\Psi(T_0)}=\inf_{\theta\in\rdue
\atop x\in\R} \| \Psi(T_0,\cdot)-(e^{i\theta_{1}}
r_{1}(\cdot-x),e^{i\theta_{2}}r_{2}(\cdot-x))\|_{\H}^2<C\eps,
$$
inequality $\Gamma_{\Psi(t)}<{\mathcal A}$ holds true by continuity for
any $t\in[T_0,T_0+\rho)$, for some small $\rho>0$, which is a contradiction by the definition
of $T_0$. Hence $T_0=T$ and, for any $T>0$, from~\eqref{firstconcl} we get
\begin{equation*}
\sup_{t\in [0,T]}\inf_{\theta\in\rdue
\atop x\in\R} \| \Psi(t,\cdot)-(e^{i\theta_{1}}
r_{1}(\cdot-x),e^{i\theta_{2}}r_{2}(\cdot-x))\|_{\H}^2<C\eps,
\end{equation*}
which is the desired property on $[0,T]$. By the arbitrariness of $T$ the assertion follows.

\bigskip
\section{Existence of a non-degenerate ground state}\label{ground}

In the following section we will show that there exists a 
non-degenerate ground state $Z$. More precisely, let us consider 
$z$ be the unique positive radial least energy solution of \eqref{eeqr} and
let $a$ be given by
\beq\label{defa}
a=(1+\beta)^{-1/2p}.
\eeq
We will prove the following result.
\bte
\label{gsnondegex}
Let $a$ be given in \eqref{defa}, then  the vector $Z=a(z,z)$ is a non-degenerate ground state  of system \eqref{ellittico} for every $p>0$, $\bt>1$ and $p\neq \bt$.
\ete
\bos
In \cite{mmp1} it is proved that for $\bt\leq 1$ every ground state of 
\eqref{ellittico} necessarily has one trivial component, that is the reason of the assumption $\bt>1$. Moreover, it can been easily seen that for $p=\bt$
the ground state $Z$ is a degenerate solution that is why we assume $p\neq \bt.$
\eos
This result will be a consequence of  the two following results.
\bte\label{esi}
Let $a$ be given in \eqref{defa}, then 
 the vector $Z=a(z,z)$ is a  ground state 
of system \eqref{ellittico} for every $p>0$, $\bt>1$.
\ete
\bte\label{deg}
Let $a$ be given in \eqref{defa}, then 
the vector $Z=a(z,z)$ is a non-degenerate  ground state 
of system \eqref{ellittico} for every $p>0$, $\bt>1$ and $p\neq \bt$.
\ete
\bos 
In \cite{fm} it is studied the global existence for the Cauchy problem 
\eqref{schr} and it is proved that the solution exists for any time 
if $p<2/n$, while it can blow up if $p\geq 2/n$. In the critical case 
$p=2/n$ it is given a bound on the $L^2$-norm of the initial data 
which guarantees the global existence of the solution (see Theorem 2).
Since Theorem \ref{esi}  shows that the test functions used in \cite{fm} 
to estimate the blow-up threshold belong to the set of ground
state solutions, as a by product, we obtain that the bound given in
\cite{fm} is the exact threshold value.
\eos
\bos
The above results have been proved for $p=1$, respectively, in~\cite{si} 
and~\cite{dw} in any dimension. Actually, the same arguments work for any $p>0$. In the following  we include the details for completeness. 
Let us notice that the same proof of Theorem \ref{esi} holds in dimension greater than one; in addition, the  arguments used in \cite{dw} hold for $p\in (0,2/n)$ for every $n\geq 1$. Thus, the vector $Z$ is a non-denerate ground state solution of \eqref{ellittico} in any dimension $n\geq 1$, our conjecture is that it is the only one if $\bt>1$.
Here our interest, is restricted to the one dimension setting so that we will
see the proof of Theorem \ref{gsnondegex} in this case.
\eos

\subsection{Proof of Theorem \ref{esi}}\label{existenceproof}

First, we recall this simple facts.

\begin{proposition}
Let us set
$$
S_{1}=\inf_{\hsob\setminus\{0\}}\frac{\|u\|^{2}_{H^1}}{\| u\|_{2p+2}^{2}},
\qquad
T_{1}=\inf_{\nehari_{1}}\Big\{\frac12\|u\|^{2}_{H^1}-\frac1{2p+2}\|u\|_{2p+2}^{2p+2}\Big\},
$$
where
$$
\nehari_{1}=\big\{u\in\hsob:\, u\neq 0,\,\,\|u\|^{2}_{H^1}=\|u\|_{2p+2}^{2p+2}\big\}.
$$
Then, the following equality holds
$$
T_{1}=\frac12\frac{p}{p+1}(S_{1})^{(p+1)/p}.
$$
\end{proposition}
\bdim
As $z$ solves the minimization problems that defines $S_{1}$ and 
$T_{1}$, using~\eqref{eeqr} we get
$$
S_{1}=\frac{\|z\|^{2}_{H^1}}{\| z\|_{2p+2}^{2}}=
\frac{\|z\|^{2}_{H^1}}{\| z\|^{2/(p+1)}}=\|z\|^{2p/(p+1)}_{H^1}=
\|z\|_{2p+2}^{2p},
$$
namely
\begin{equation}
\label{euqaS}
\|z\|^{2}_{H^1}=S_{1}^{(p+1)/p}
\,\,\quad
\text{and} 
\,\,\quad
\|z\|_{2p+2}=S_{1}^{1/2p}.
\end{equation}
Using these equalities in the definition of $T_{1}$ permits to conclude the proof.
\edim

\vskip2pt
\noindent
Define now the sets
\begin{align*}
{\mathcal N}_{0} &=\left\{U\in\H: U\neq (0,0),\,\, 
\|U\|^{2}_{\H}=\|U\|_{2p+2}^{2p+2}+2\bt\|u_{1}u_{2}\|_{p+1}^{p+1}\right\}, \\
{\mathcal N} &=\big\{
U\in\H: u_{i}\neq 0,\,\,
 \|u_{i}\|^{2}_{H^1}=\|u_{i}\|_{2p+2}^{2p+2}+
\bt\| u_{1}u_{2}\|^{p+1}_{p+1},\,\,i=1,2\big\}.
\end{align*}
Moreover, if $\H_{r}$ is the set of radial function of $\H$, we introduce the numbers
\begin{equation}
\label{definf}
A_{0}=\inf_{U\in \nehari_{0}} \I(U),\qquad A=\inf_{U\in \nehari} \I(U),
\qquad A_{r}=\inf_{U\in \nehari\cap \H_{r}} \I(U),
\end{equation}
where
$$
\I(U)=\frac{1}{2}\|U\|^{2}_{\H}-\frac{1}{2p+2}\|U\|_{2p+2}^{2p+2}
-\frac{1}{p+1}\bt\|u_{1}u_{2}\|_{p+1}^{p+1}.
$$
Let $a$ be a positive number. Writing down the equations that define 
$\nehari$ and recalling that $z$ satisfies \eqref{eeqr} it is easy
to see that $a(z,z)\in \nehari$ if $a$ satisfies \eqref{defa}.
\\
\noindent
Concerning the infimum problems $A_0,A,A_r$, in \cite{si} the  following
result is proved for $p=1$; actually the same proof holds for any $p$
satisfying \eqref{pzero}, we include some details.
\begin{proposition}\label{disinf}
Let $a$ satisfies \eqref{defa}. Then the following inequalities hold
\beq\label{disA}
0<A_{0}\leq A\leq A_{r}\leq \frac{p}{p+1}a^{2}S_{1}^{(p+1)/p},
\eeq
where the values $A_{0}$ and  $A_{r}$ are defined in~\eqref{definf}.
\end{proposition}

\bdim
First note that, taken any $U=(u_1,u_2)\in \nehari_{0} $, the value $\I(U)$ is equal to
\begin{equation}\label{disI}
\I(U)=\frac12\Big(\frac{p}{p+1}\Big)\big[\|U\|_{2p+2}^{2p+2}+2\beta
\|u_{1}u_{2}\|_{p+1}^{p+1}\big]=\frac12\Big(\frac{p}{p+1}\Big)\|U\|^{2}_{\H}.
\end{equation}
Moreover, since $a(z,z)\in \nehari$ and has radial components, recalling~\eqref{euqaS} we get
\beq\label{disar}
A_{r} \leq \I(az,az) = \frac12\Big(\frac{p}{p+1}\Big)
\|(az,az)\|^{2}_{H^1}=\Big(\frac{p}{p+1}\Big)a^{2}\|z\|^{2}_{H^1}  
=\Big(\frac{p}{p+1}\Big)a^{2} S_{1}^{(p+1)/p},
\eeq
which is the last inequality on the right-hand side in \eqref{disA}. 
It just remains to show that $A_{0}>0$. To this aim, 
take $U\in \nehari_{0}$  and observe that H\"older and Sobolev 
inequalities imply that there exist positive constants $C_0,C_1$ such that
$$
\|U\|^{2}_{\H}=\|U\|_{2p+2}^{2p+2}+2\beta\|u_{1}u_{2}\|_{p+1}^{p+1}\leq 
C_0\|U\|_{2p+2}^{2p+2}\leq C_1\|U\|^{2p+2}_{\H}
$$
so that the norm $\|U\|_{\H}$ remains uniformly away from zero. Hence, 
recalling formula~\eqref{disI}, we conclude the proof.
\edim

\vskip4pt

We are now ready to complete the proof of Theorem \ref{esi}.
\\
\noindent{\bf Proof of Theorem~\ref{esi} concluded.}
We will obtain Theorem \ref{esi} by showing that the infimum $A$ equals $A_{r}$ and it is achieved at the couple $a(z,z)$, which is thus a 
ground state solution of~\eqref{ellittico}.
\\
First, let $(U_{m})=(u_{m,1},u_{m,2})\subset {\mathcal N}$ be a minimizing sequence for $A$, namely
$\I(U_m)=A+o(1)$ as $m\to\infty$. Let us set $y_{m,i}=\|u_{m,i}\|_{2p+2}^{2}$ for any 
$m\in\N$ and $i=1,2$. Hence, by the definition of $S_1$ and H\"older
inequality, it follows that, for all $m\in\N$,
\begin{equation}
\label{disy}
S_{1}y_{m,1}\leq \|\unu\|^{2}_{H^1}=\|\unu\|_{2p+2}^{2p+2}+\beta 
\|\unu\und\|_{p+1}^{p+1}\leq \ynu^{p+1}+\beta \ynu^{(p+1)/2}
\ynd^{(p+1)/2},
\end{equation}
for all $m\in\N$. Of course, for all $m\in\N$, the analogous inequality holds
\begin{equation}
\label{disy-2}
S_{1}y_{m,2}\leq \|\und\|^{2}_{H^1}=\|\und\|_{2p+2}^{2p+2}+\beta 
\|\unu\und\|_{p+1}^{p+1}\leq \ynd^{p+1}+\beta \ynu^{(p+1)/2}
\ynd^{(p+1)/2}.
\end{equation}
Furthermore, taking into account formula~\eqref{disI}, by addition of the first 
inequalities in~\eqref{disy} and~\eqref{disy-2} one obtains
\beq\label{disfina}
S_{1}(\ynu+\ynd)\leq 2\frac{p+1}{p}\I (U_{n})=2\frac{p+1}{p} A+o(1),\quad\text{as $m\to\infty$}.
\eeq
By combining this inequality with Proposition~\ref{disinf} gives
$$
S_{1}(\ynu+\ynd)\leq 2a^{2}S_{1}^{(p+1)/p}+o(1),\quad\text{as $m\to\infty$}.
$$
Hence, defining $z_{m,i}=y_{m,i}/S_{1}^{1/p}$, we derive
$z_{m,1}+z_{m,2}\leq 2a^{2}+o(1),$ as $m$ tends to infinity.
Also, by dividing \eqref{disy} by $S_{1}\ynu$ and~\eqref{disy-2} by $S_{1}\ynd$ and using 
$S_{1}=S_{1}^{(p-1)/2p}S_{1}^{(p+1)/2p}$ we obtain
that, as $m\to\infty$,  $(z_{m,1},z_{m,2})$ satisfies the following system of inequalities
$$
\begin{cases}
z_{m,1}+z_{m,2}\leq 2a^{2}+o(1),
\medskip\\
z^{p}_{m,1}+\bt z_{m,1}^{(p-1)/2}z_{m,2}^{(p+1)/2}  \geq 1,
\medskip\\
z^{p}_{m,2}+\bt z_{m,1}^{(p+1)/2}z_{m,2}^{(p-1)/2} \geq 1.
\end{cases}
$$
Taking into account \eqref{defa} we are lead to the 
study of the associated algebraic system of inequalities
\begin{equation}\label{algebinequalit}
\begin{cases}
x+y \leq 2a^{2},
\medskip\\
x^{p}+\bt x^{(p-1)/2}y^{(p+1)/2} \geq (1+\beta)a^{2p},
\medskip\\
y^{p}+\bt x^{(p+1)/2}y^{(p-1)/2} \geq (1+\beta)a^{2p},
\end{cases}
\end{equation}
for which we refer to Figure 1.

\noindent
Then, for $\beta >1$ and any $i=1,2$, the sequence $(z_{m,i})$ remains bounded away from zero and it has to be $z_{m,1}\to a^{2}$ and $z_{m,2}\to a^{2}$ as $m\to\infty$, so that looking at the first (in)equality of~\eqref{algebinequalit} with $x=y$ (by figure 1) yields $x=y=a^2$), so that 
$\ynu\to a^{2}S_{1}^{1/p},$ and 
$\ynd\to a^{2}S_{1}^{1/p},$ as $m$ diverges.
Whence, passing to the limit in formula~\eqref{disfina}, in light of Proposition \ref{disinf}
we obtain
$$
2S_{1}^{(p+1)/p}a^{2}\leq 2\frac{p+1}{p} A\leq 2a^{2}S_{1}^{(p+1)/p}
$$
so that, \eqref{disar}, gives
$$
A\leq A_{r} \leq \I(az,az) \leq \Big(\frac{p}{p+1}\Big)a^{2}
\left(S_{1}\right)^{(p+1)/p}=A,
$$
which gives $A=A_{r}=\I(az,az)$, concluding the proof.
\edim

\subsection{Proof of Theorem \ref{deg}}

According to Section~\ref{existenceproof}, let us consider $Z=a(z,z)$ the
particular ground state solution of \eqref{ellittico}, with $a$ given
in \eqref{defa}; we will now show the non-degeneracy property of $Z$.
First, notice that the linearized system \eqref{linsyste} can be obtained using the operator $L_{+}$ acting on $Z$, and by the explicit expression of $Z$ we
get 
\begin{equation*}
L_{+}=
\left(
\begin{array}{cc}
-\dfrac12\pa_{xx} +1 & 0
\\\\
0 & -\dfrac12\pa_{xx} +1
\end{array}
\right)-\left(
\begin{array}{cc}
\dfrac{p(2+\bt)+1}{1+\bt}z^{2p} 
& \dfrac{\bt(p+1)}{1+\bt}z^{2p} \\\\
\dfrac{\bt(p+1)}{1+\bt}z^{2p} & \dfrac{p(2+\bt) +1}{1+\bt}z^{2p}
\end{array}
\right).
\end{equation*}
In accordance with Section \ref{spectral}, we denote with $H_{F}(Z)$
the second matrix on the right hand side. The quadratic form related to $H_{F}(Z)$ can be diagonalized by an orthonormal change 
of coordinates, introducing
\begin{equation}\label{coordinatech}
w_{1}=\dfrac{\sqrt{2}}{2}(\phi_{1}+\phi_{2}),\quad
w_{2}=\dfrac{\sqrt{2}}{2}(\phi_{1}-\phi_{2}).
\end{equation}
Since we have
\bdm
\hbox{Tr}(H_{F}(Z))=2 \dfrac{(2+\bt)p+1}{1+\bt}
=(2p+1)+\dfrac{2p+1-\bt}{1+\bt},\qquad
\hbox{Det}(H_{F}(Z))=\dfrac{(2p+1)(2p+1-\bt)}{1+\bt},
\edm
it follows that its eigenvalues are
\beq\label{auto}
\la_{1}=2p+1,\qquad \la_{2}=\dfrac{2p+1-\bt}{1+\bt}\in
(-1,2p+1)
\eeq
so the linear elliptic system $L_{+}\Phi=0$ decouples and reduces to
\beq\label{diagonale}
\begin{cases}
-\frac12\pa_{xx} w_{1}+w_{1}=(2p+1)z^{2p}(x) w_{1}, & \text{in $\R$}\\ 
-\frac12\pa_{xx} w_{2}+w_{2}=\dfrac{2p+1-\bt}{1+\bt}z^{2p}(x)  w_{2}, & \text{in $\R$}.
\end{cases}
\eeq
Taking into account that the weight $z$ is exponentially decaying, the spectrum
of the linear self-adjoint operator $-\frac{1}{2}\partial_{xx}+{\rm Id}-\mu z^{2p}$ is discrete.
Furthermore, from \cite[(a) and (b) of Proposition 2.8]{weinsteinMS} with proofs for $n=1$
in~\cite[Appendix A]{weinsteinMS}, we learn that the eigenvalues of 
\beq\label{eigen}
-\dfrac12\pa_{xx} w+w-\mu z^{2p}(x) w=0\qquad \text{in $\R$},
\eeq
are given by $\mu_{1}=1,$  $\mu_{2}=2p+1,$  $\mu_{3}>2p+1,$
and, denoting by $V_{\mu_i}$ the eigenspace corresponding to the eigenvalue $\mu_i$, we have $ V_{\mu_1}={\rm span}\big\{z\big\},$
$ V_{\mu_2}={\rm span}\big\{\partial_x z\big\}.$
Therefore, from the first equation of~\eqref{diagonale} we deduce
$w_1\in {\rm span}\big\{\partial_x z\big\}.$ From \eqref{auto} we also
deduce, from the second equation of~\eqref{diagonale}, that
$w_2=0.$
In turn, by the orthonormal change of coordinates~\eqref{coordinatech} we obtain  $\phi_1=\phi_2=c\partial_x z$, for some coefficient $c\in\R$. Whence 
$\hbox{Ker}(L_{+})=\langle\pa_{x}Z_{\bt}\rangle$, which concludes the proof.
\edim

\end{document}